\documentclass{amsart}
\usepackage{amsmath,amssymb}
\usepackage[dvips]{graphics}
\usepackage[all]{xy}
\newtheorem{Theorem}{Theorem}[section]

\newtheorem{Proposition}[Theorem]{Proposition}

\newtheorem{Example}[Theorem]{Example}

\newtheorem{Remark}[Theorem]{Remark}

\makeatletter
\@addtoreset{figure}{section}
\def\@thmcountersep{-}
\makeatother

\numberwithin{equation}{section}

\begin{document}

\title{Braid presentation of spatial graphs}

\author{Ken Kanno}
\address{Graduate School of Education, Waseda University, Nishi-Waseda 1-6-1, Shinjuku-ku, Tokyo, 169-8050, Japan}
\email{kanno@suou.waseda.jp}

\author{Kouki Taniyama}
\address{Department of Mathematics, School of Education, Waseda University, Nishi-Waseda 1-6-1, Shinjuku-ku, Tokyo, 169-8050, Japan}
\email{taniyama@waseda.jp}

\subjclass[2000]{Primary 57M25; Secondary 57M15}

\date{}

\dedicatory{}

\keywords{spatial graph, braid presentation}

\begin{abstract}
We define braid presentation of edge-oriented spatial graphs as a natural generalization of braid presentation of oriented links. We show that every spatial graph has a braid presentation. For an oriented link it is known that the braid index is equal to the minimal number of Seifert circles. We show that an analogy does not hold for spatial graphs.
\end{abstract}

\maketitle

\section{Introduction} 

Throughout this paper we work in the piecewise linear category. Let $G$ be a finite edge-oriented graph. Namely $G$ consists of finite vertices and finite edges, and each edge has a fixed orientation. Edge-oriented graph is called digraph in graph theory. We consider a graph as a topological space in a usual way. Let ${\mathbb S}^3$ be the unit 3-sphere in the $xyzw$-space ${\mathbb R}^4$ centered at the origin of ${\mathbb R}^4$. An embedding of $G$ into ${\mathbb S}^3$ is called a {\it spatial embedding} of $G$. Then the image is also called a spatial embedding or a {\it spatial graph}. 
Let $A$ (resp. $C$) be the intersection of ${\mathbb S}^3$ and the $zw$-plane (resp. $xy$-plane). Then the union $A\cup C$ is a Hopf link in ${\mathbb S}^3$. We call $A$ the {\it axis} and $C$ the {\it core}. 
Let $\pi:{\mathbb S}^3-A\to C$ be a natural projection defined by 
$\displaystyle{\pi(x,y,z,w)=(\frac{x}{\sqrt{x^2+y^2}},\frac{y}{\sqrt{x^2+y^2}},0,0)}$. 
We give a counter-clockwise orientation to $C$ on $xy$-plane and fix it. We say that a continuous map $\varphi:G\to C$ is {\it locally orientation preserving} if for any edge $e$ of $G$ and any point $p$ on $e$ there is a neighbourhood $U$ of $p$ in $e$ such that the restriction map of $\varphi$ to $U$ is an orientation preserving embedding. Let $f:G\to {\mathbb S}^3$ be a spatial embedding. We say that $f$ or its image $f(G)$ is a {\it braid presentation} if $f(G)$ is disjoint from $A$ and the composition map $\pi\circ f':G\to C$ is locally orientation preserving where $f':G\to {\mathbb S}^3-A$ is the map defined by $f'(p)=f(p)$ for any $p$ in $G$. Note that this generalizes the braid presentation defined for $\theta_m$-curve in \cite{S-T}. 
The following theorem shows that every edge-oriented spatial graph has a braid presentation up to ambient isotopy of ${\mathbb S}^3$. This generalizes Alexander's theorem that every oriented link can be expressed by a closed braid \cite{Alexander} and that proved for $\theta_m$-curve in \cite{S-T}. 

\vskip 3mm

\begin{Theorem}\label{main-theorem1}
Let $G$ be a finite edge-oriented graph and $f:G\to {\mathbb S}^3$ a spatial embedding. Then there is a braid presentation $g:G\to {\mathbb S}^3$ that is ambient isotopic to $f$ in ${\mathbb S}^3$.
\end{Theorem}

\vskip 3mm

In \cite{Yamada} it is shown that the minimal number of Seifert circles of an oriented link is equal to the braid index of the link. We consider an analogy to spatial graphs. Let $g:G\to {\mathbb S}^3$ be a braid presentation. Let $P_p=\pi^{-1}(p)$ for $p\in C$. Let $\tilde{b}(g)=\tilde{b}(g(G))$ be the maximum of $|g(G)\cap P_p|$ where $|X|$ denotes the cardinality of the set $X$ and $p$ varies over all points in $C$. Let $f:G\to {\mathbb S}^3$ a spatial embedding. Let $b(f)=b(f(G))$ be the minimum of $\tilde{b}(g)$ where $g$ varies over all braid presentation that is ambient isotopic to $f$.
We call $b(f)=b(f(G))$ the {\it braid index} of $f$ or $f(G)$. 

\begin{figure}[htbp]
      \begin{center}
\scalebox{0.4}{\includegraphics*{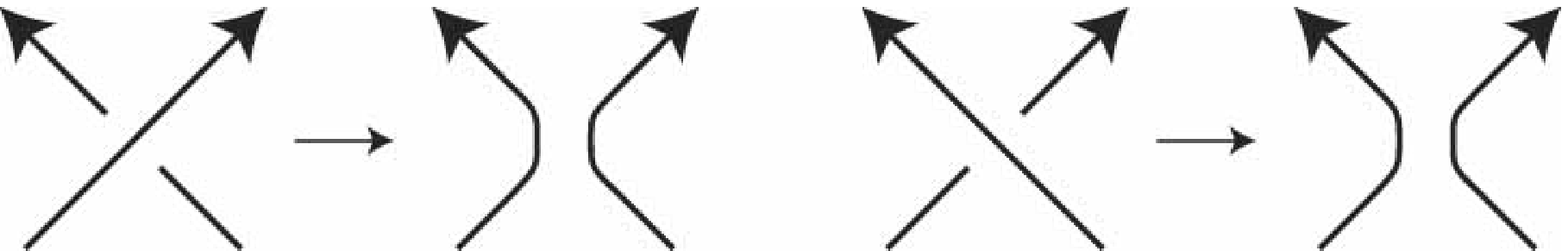}}
      \end{center}
   \caption{}
  \label{smoothing}
\end{figure}

Let ${\mathbb S}^2$ be the intersection of ${\mathbb S}^3$ and the $xyz$-space. Then any spatial embedding $f:G\to {\mathbb S}^3$ has a diagram on ${\mathbb S}^2$ up to ambient isotopy of ${\mathbb S}^3$. Let $D$ be a diagram of $f$ on ${\mathbb S}^2$. Let $S(D)$ be a plane graph in ${\mathbb S}^2$ obtained from $D$ by smoothing every crossing of $D$. Here smoothing respect the orientations of the edges. See Figure \ref{smoothing}.

%

Let $\mu(X)$ be the number of connected components of a space $X$. Let $s(f)=s(f(G))$ be the minimum of $\mu(S(D))$ where $D$ varies over all diagrams of $f$ up to ambient isotopy of ${\mathbb S}^3$. We call $s(f)=s(f(G))$ the {\it smoothing index} of $f$ or $f(G)$. Note that our $s(f)$ is different from $s(G)$ defined for $\theta_m$-curve in \cite{S-T}. 
We will show in Proposition \ref{proposition1} that for any natural number $n$ there is a spatial embedding $f$ of $G$ with $b(f)\geq n$ unless $G$ contains no cycles as an unoriented graph. In contrast we will show in the next theorem that $s(f,G)=s(g,G)$ for any two spatial embeddings $f$ and $g$ of $G$ unless $G$ satisfies certain conditions. 
By ${\rm indeg}(v,G)={\rm indeg}(v)$ (resp. ${\rm outdeg}(v,G)={\rm outdeg}(v)$) we denote the number of the edges whose head (resp. tail) is the vertex $v$ of $G$. Then ${\rm deg}(v,G)={\rm deg}(v)={\rm indeg}(v,G)+{\rm outdeg}(v,G)$ is called the {\it degree} of $v$ in $G$. 
We say that an edge-oriented graph $G$ is {\it circulating} if ${\rm indeg}(v,G)={\rm outdeg}(v,G)$ for any vertex $v$ of $G$. Namely each component of a circulating graph is Eulerian. Let $\chi(X)$ be the Euler characteristic of a space $X$. 
Then we have the following theorem.

\vskip 3mm

\begin{Theorem}\label{main-theorem2}
Let $G$ be a finite edge-oriented graph without isolated vertices.

(1) Suppose that $G$ is not circulating. Then for any spatial embedding $f:G\to {\mathbb S}^3$, $s(f)={\rm max}\{1,\chi(G)\}$.

(2) Suppose that $G$ is circulating. Then for any natural number $n$ there is a spatial embedding $f:G\to {\mathbb S}^3$ such that $s(f)\geq n$.

\end{Theorem}

\vskip 3mm

\begin{Remark}\label{remark}
{\rm
Another choice for the smoothing index as a generalization of the number of Seifert circles of an oriented link is the use of the first Betti number instead of the number of connected components. Let $s'(f)$ be the minimum of $\beta_1(S(D))$ among all diagrams $D$ of $f$. By the Euler-Poincar\'{e} formula we have $\beta_1(S(D))=\mu(S(D))-\chi(S(D))$. Since smoothing does not change the Euler characteristic we have $\chi(S(D))=\chi(G)$. Then we have $s'(f)=s(f)-\chi(G)$. Thus we have that $s'(f)$ is determined by $s(f)$ after all.
}
\end{Remark}

\section{Proof of Theorem \ref{main-theorem1}} 
The following proof is a natural extension of a proof of Alexander's theorem by Cromwell \cite{Cromwell} using {\it rectangular diagram} of oriented links that appears in \cite{Brunn} \cite{Cromwell} \cite{M-M} etc. 
In this section we regard ${\mathbb S}^2$ as a one-point compactification of the $xy$-plane. Thus we may suppose that all diagrams are on the $xy$-plane. 
In the following we sometimes do not distinguish an abstract vertex or edge from its image in ${\mathbb S}^3$ or on ${\mathbb S}^2$.

\vskip 3mm
\noindent{\bf Proof of Theorem \ref{main-theorem1}.}
Let $D$ be a diagram of the spatial embedding $f:G\to {\mathbb S}^3$. We will deform $D$ step by step so that it is still a diagram of $f$ up to ambient isotopy in ${\mathbb S}^3$ as follows. First we move $D$ if necessary so that $D$ is left to the $y$-axis. Namely $D$ is contained in the region of the $xy$-plane defined by $x<0$. 
By a local deformation near each vertex we may suppose that all edges go down with respect to the $y$-coordinate in each small neighbourhood of a vertex of $G$. Then we further deform $D$ so that it satisfies the following conditions.

(1) $D$ is a union of finitely many line segments.

(2) Each vertex $v$ has a small disk neighbourhood $N_v$ such that the diagram $D$ on $N_v$ is a union of ${\rm indeg}(v,G)+{\rm outdeg}(v,G)$ line segments each of which has $v$ as one of its end points, and each of which goes down with respect to the $y$-coordinate.

\begin{figure}[htbp]
      \begin{center}
\scalebox{0.4}{\includegraphics*{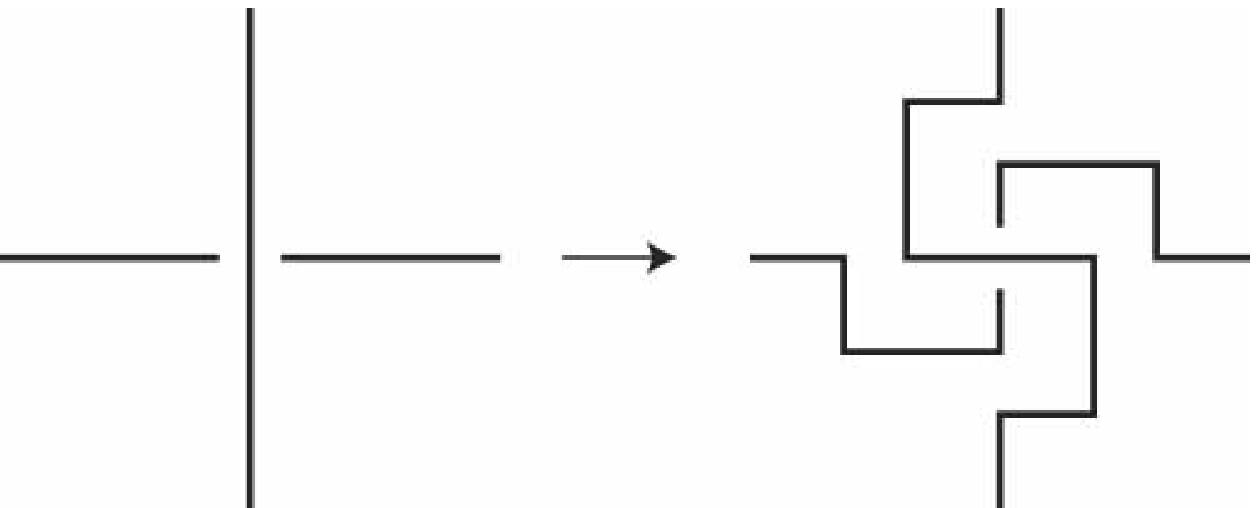}}
      \end{center}
   \caption{}
  \label{manji}
\end{figure}

\begin{figure}[htbp]
      \begin{center}
\scalebox{0.4}{\includegraphics*{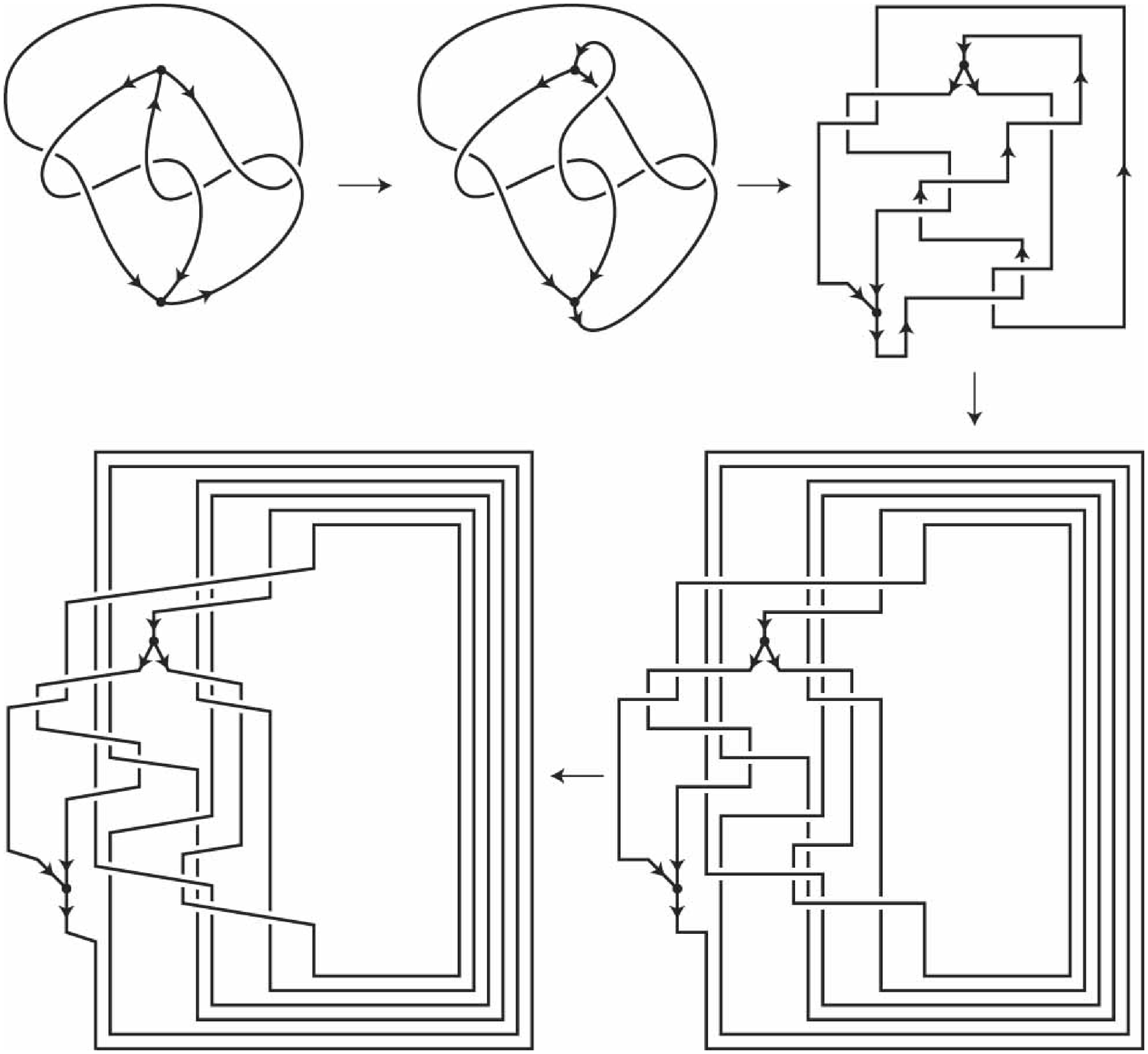}}
      \end{center}
   \caption{}
  \label{braiding}
\end{figure}

(3) A line segment that is not contained in any $N_v$ is parallel to $x$-axis or $y$-axis.
Then we have that each crossing of $D$ is a crossing between a horizontal line segment (parallel to $x$-axis) and a vertical line segment (parallel to $y$-axis). Then by a local deformation as illustrated in Figure \ref{manji} we may further assume that a horizontal line segment is over a vertical line segment at each crossing. Note that the disk neighbourhood $N_v$ can be taken to be arbitrarily small. Then by a slight deformation we have that a straight line that contains a vertical line segment going up with respect to the $y$-coordinate contains no other vertical line segments and is disjoint from any $N_v$.
Let $s_1,\cdots,s_n$ be the vertical line segments that go up with respect to the $y$-coordinate. We may suppose without loss of generality that the $x$-coordinate of $s_i$ is less than that of $s_j$ if $i<j$. Let $R_1,\cdots,R_n$ be sufficiently large upright rectangles with $R_1\supset\cdots\supset R_n$ and $s_i\subset \partial R_i$ for each $i$. We replace each $s_i$ by $\partial R_i-{\rm int}(s_i)$ that crosses under the horizontal line segments at every crossings. Finally we tilt the horizontal line segments other than $\partial R_i$ slightly so that they go down with respect to the $y$-coordinate. Then we finally have a diagram $D$ of $f$ that totally turns around the origin of the $xy$-plane. Then $D$ represents a braid presentation. See for example Figure \ref{braiding}. This completes the proof. $\Box$

%

%

%

%
\section{Proof of Theorem \ref{main-theorem2}} 
A vertex $v$ of an edge-oriented graph $G$ is called a {\it source} (resp. {\it sink}) of $G$ if ${\rm indeg}(v,G)=0$ (resp. ${\rm outdeg}(v,G)=0$).

\begin{Proposition}\label{proposition1}
Let $G$ be a finite edge-oriented graph. Suppose that $G$ contains a cycle as an unoriented graph. 
Then for any natural number $n$ there is a spatial embedding $f:G\to {\mathbb S}^3$ such that $b(f)\geq n$.
\end{Proposition}

\vskip 3mm
\noindent{\bf Proof.} Let $\gamma$ be a cycle of $G$. Note that $\gamma$ may not be an oriented cycle as an edge-oriented subgraph of $G$. Let $\alpha$ be the number of the sources of $\gamma$. Let $f:G\to {\mathbb S}^3$ be a spatial embedding of $G$ such that the bridge index ${\rm bridge}(f(\gamma))$ of the knot $f(\gamma)$ is greater than or equal to $n+\alpha$. By the definition we have $b(f(\gamma))\leq b(f(G))$. We may suppose that $f(\gamma)$ is a braid presentation with $\tilde{b}(f(\gamma))=b(f(\gamma))$. Then $f(\gamma)$ has at most $2(\tilde{b}(f(\gamma))+\alpha)$ critical points with respect to the $y$-coordinate. Therefore we have ${\rm bridge}(f(\gamma))\leq\tilde{b}(f(\gamma))+\alpha$. Therefore we have $n+\alpha\leq{\rm bridge}(f(\gamma))\leq\tilde{b}(f(\gamma))+\alpha=b(f(\gamma))+\alpha$. Thus we have $n\leq b(f(\gamma))\leq b(f(G))$ as desired. $\Box$

\vskip 3mm

\begin{Proposition}\label{proposition2}
Let $G$ be a finite edge-oriented graph and $f:G\to {\mathbb S}^3$ a spatial embedding.Then $s(f)\geq\chi(G)$.
\end{Proposition}

\vskip 3mm
\noindent{\bf Proof.} Let $D$ be a diagram of $f$. It is sufficient to show that $\mu(S(D))\geq\chi(G)$. 
Since smoothing does not change the Euler characteristic we have that $\chi(G)=\chi(S(D))$. 
Let $\beta_1(X)$ be the first Betti number of a space $X$. 
By the Euler-Poincar\'{e} formula we have $\chi(S(D))=\mu(S(D))-\beta_1(S(D))$. Therefore we have $\chi(S(D))\leq\mu(S(D))$. Thus we have $\mu(S(D))\geq\chi(G)$. $\Box$

\vskip 3mm

\noindent{\bf Proof of Theorem \ref{main-theorem2} (1).} First we show that $s(f)\geq{\rm max}\{1,\chi(G)\}$. Let $D$ be a diagram of $f$. By the definition we have $\mu(S(D))\geq1$. By Proposition \ref{proposition2} we have $\mu(S(D))\geq\chi(G)$. Thus we have $\mu(S(D))\geq{\rm max}\{1,\chi(G)\}$ holds for any diagram $D$ of $f$. This implies $s(f)\geq{\rm max}\{1,\chi(G)\}$.


Next we show that $f(G)$ has a diagram $D$ with $\mu(S(D))={\rm max}\{1,\chi(G)\}$. Let $H$ be the maximal subgraph of $G$ that has no vertices of degree less than 2. If $H$ is not an empty graph and $H$ is not circulating then we set $G'=H$. Suppose that $H$ is an empty graph or a circulating graph. Suppose that there is a component $I$ of $H$ that is not a component of $G$. Let $J$ be the component of $G$ containing $I$. Let $e$ be an edge of $J$ that is not an edge of $I$ but incident to a vertex of $I$. Let $G'$ be the minimal subgraph of $G$ that contains $H$ and $e$. Suppose that every component of $H$ is also a component of $G$. Let $e$ be an edge of $G$ that is not an edge of $H$. Let $G'$ be the minimal subgraph of $G$ that contains $H$ and $e$. Note that in any case we have $\chi(G')\leq1$.
Let $f'$ be the restriction map of the spatial embedding $f:G\to {\mathbb S}^3$ to $G'$. We will construct a diagram $D'$ of $f'$ with $\mu(S(D'))={\rm max}\{1,\chi(G')\}=1$. Namely we will construct $D'$ such that $S(D')$ is connected. We start from a diagram $D'$ of $f'$ and deform it step by step and finally have $D'$ with $S(D')$ connected. 
By Theorem \ref{main-theorem1} we may suppose that $f'$ is a braid presentation. By deforming the braid presentation if necessary we have that $f'$ has a diagram $D'$ on the $xy$-plane with the following properties.

(1) There exists a rectangle $B=[-3,-1]\times[-2,2]$ in the $xy$-plane such that at every point on $D'$ in $B$ the edge-orientation goes down with respect to the $y$-coordinate.

(2) Outside of $B$ the diagram $D'$ consists of some parallel arcs turning around the origin of the $xy$-plane.

\begin{figure}[htbp]
      \begin{center}
\scalebox{0.6}{\includegraphics*{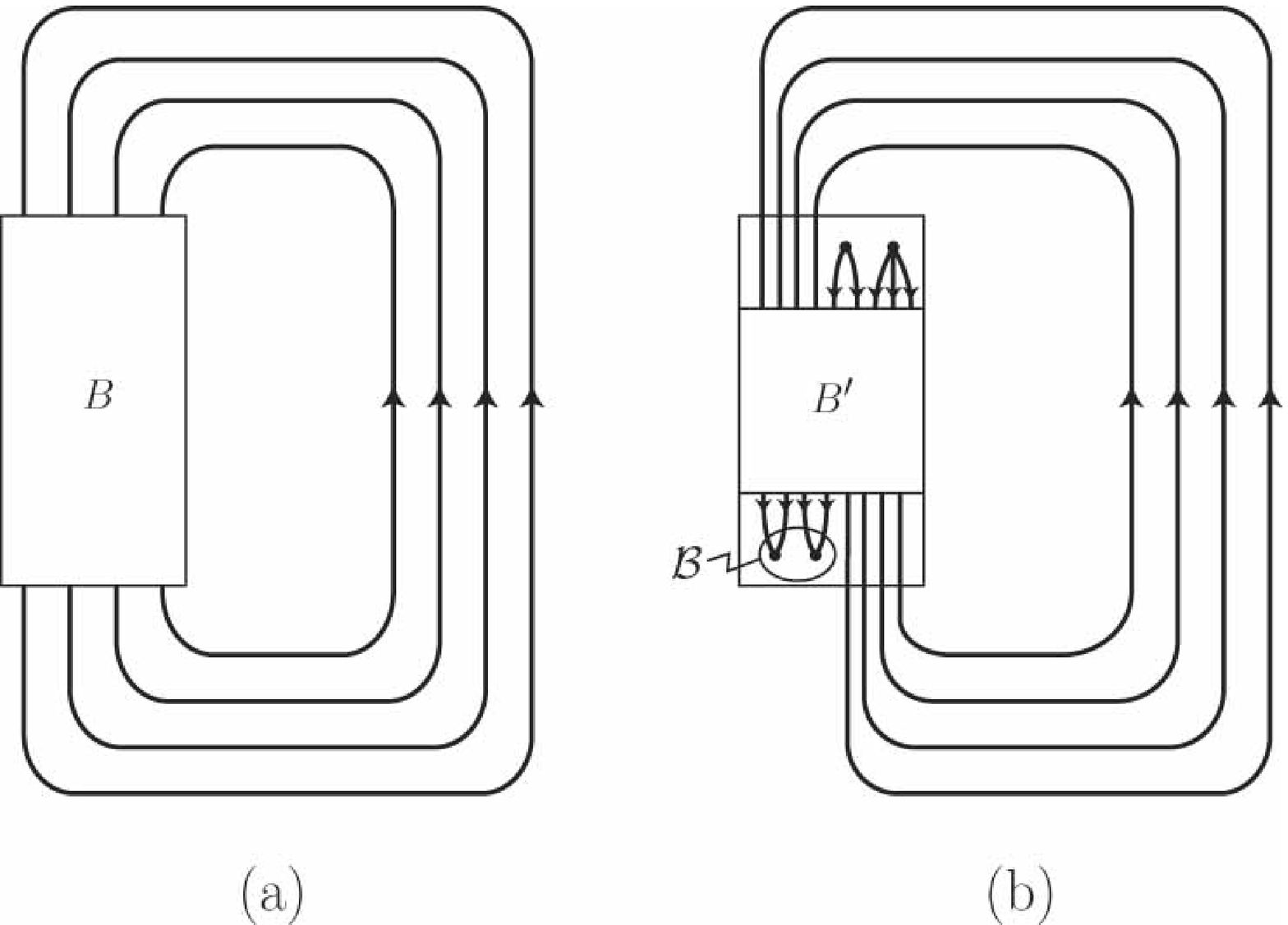}}
      \end{center}
   \caption{}
  \label{proof1}
\end{figure}

See for example Figure \ref{proof1} (a). In the following deformations we always keep the condition that everything goes down with respect to the $y$-coordinate inside $B$.

\begin{figure}[htbp]
      \begin{center}
\scalebox{0.4}{\includegraphics*{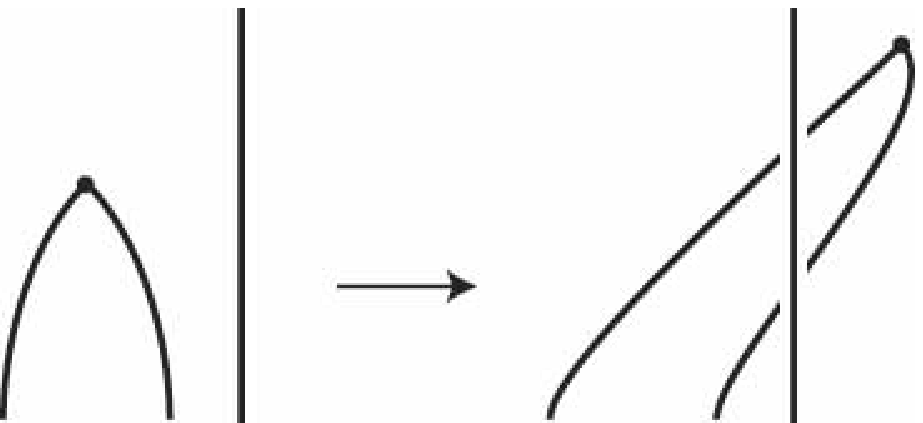}}
      \end{center}
   \caption{}
  \label{proof1-2}
\end{figure}

Suppose that $G'$ has some sources and/or sinks. Then by pulling up the sources and moving them to the right as illustrated in Figure \ref{proof1-2} and pulling down the sinks and moving them to the left we have that all sources are in $[-3,-1]\times[1,2]$ and all sinks are in $[-3,-1]\times[-2,-1]$ and all parallel arcs go left to the sources and go right to the sinks outside of $B'=[-3,-1]\times[-1,1]$. See for example Figure \ref{proof1} (b).

%

%

%

Then we deform the diagram $D'$ in $B'$ such that the following conditions hold.

(1) If a vertex $v$ of $G'$ satisfies $1\leq{\rm indeg}(v,G')\leq{\rm outdeg}(v,G')$ then it is rightmost in $B'$.

(2) If a vertex $v$ of $G'$ satisfies ${\rm indeg}(v,G')>{\rm outdeg}(v,G')\geq1$ then it is leftmost in $B'$.

\begin{figure}[htbp]
      \begin{center}
\scalebox{0.7}{\includegraphics*{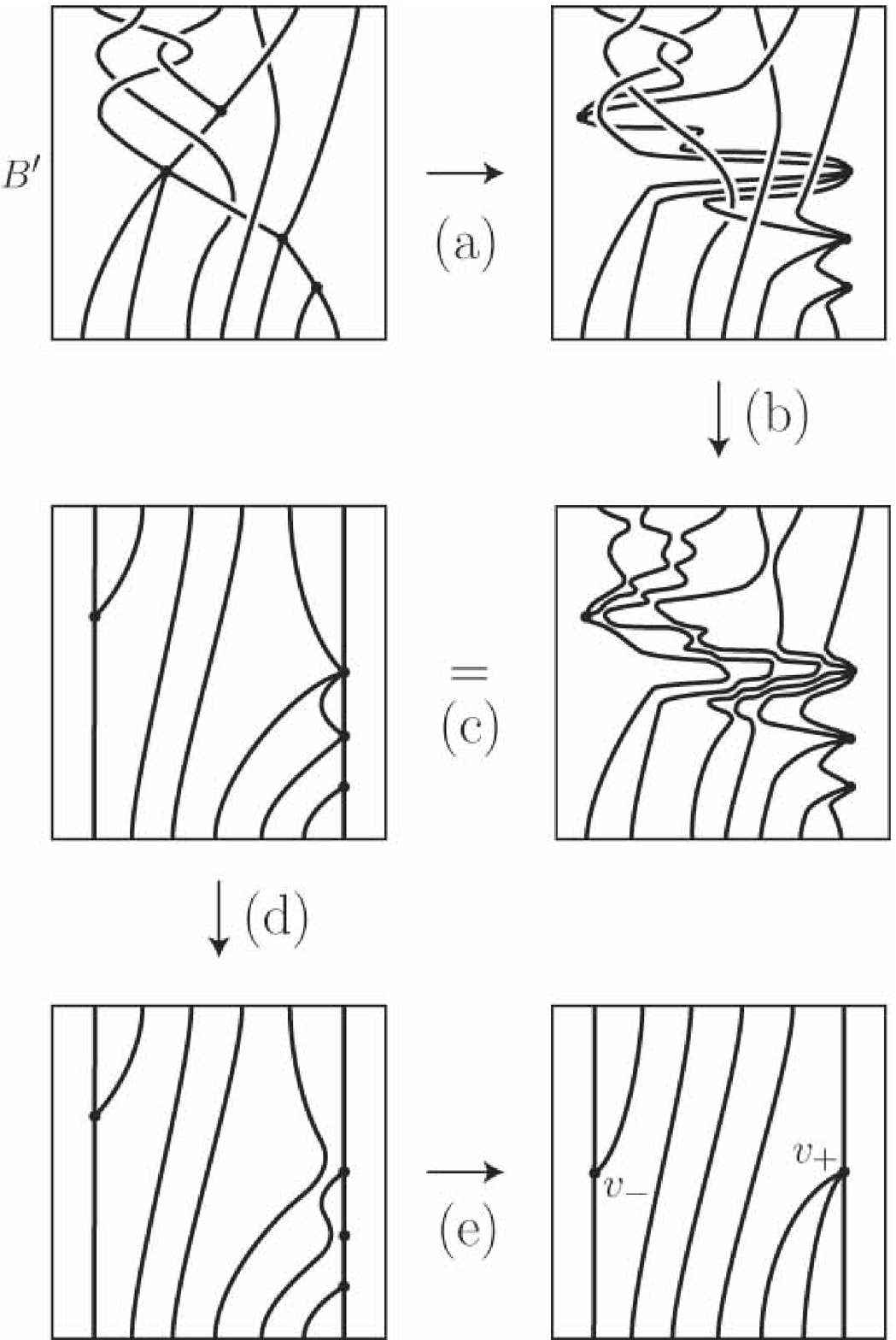}}
      \end{center}
   \caption{}
  \label{proof2}
\end{figure}

\begin{figure}[htbp]
      \begin{center}
\scalebox{0.6}{\includegraphics*{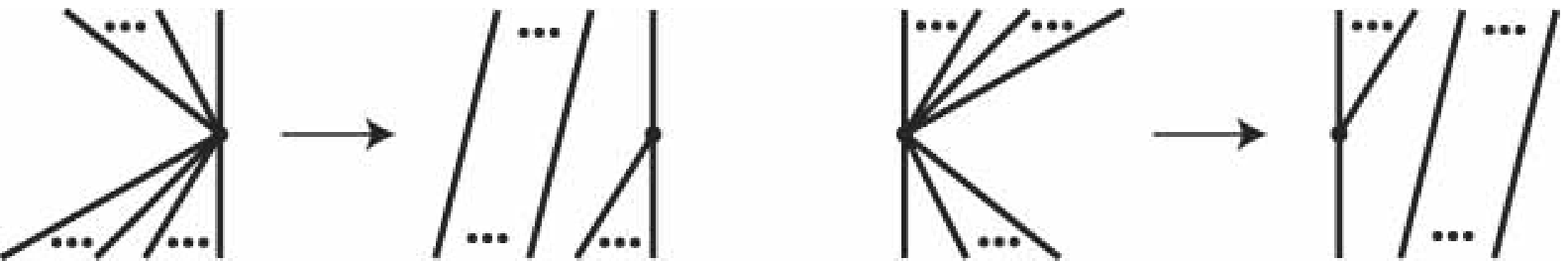}}
      \end{center}
   \caption{}
  \label{proof3}
\end{figure}

See for example Figure \ref{proof2} (a). 
Now we perform smoothing for the crossings in $B'$ and obtain $S(D')$ on $B'$. See for example Figure \ref{proof2} (b) and (c). 
We do not deform $D'$ inside $B'$ any more. We will deform $D'$ only outside of $B'$. However to make the situation simple we further perform the following replacement of $S(D')$ on $B'$. For each vertex $v$ with $2\leq{\rm indeg}(v,G')\leq{\rm outdeg}(v,G')$ we replace a neighbourhood of it on $B'$ by ${\rm indeg}(v,G')-1$ parallel arcs and a vertex $u$ with ${\rm indeg}(u)=1$ and ${\rm outdeg}(u)={\rm outdeg}(v,G')-{\rm indeg}(v,G')+1$. Similarly for each vertex $v$ with ${\rm indeg}(v,G')>{\rm outdeg}(v,G')\geq2$ we replace a neighbourhood of it on $B'$ by ${\rm outdeg}(v,G')-1$ parallel arcs and a vertex $u$ with ${\rm outdeg}(u)=1$ and ${\rm indeg}(u)={\rm indeg}(v,G')-{\rm outdeg}(v,G')+1$. See Figure \ref{proof3} and Figure \ref{proof2} (d). Then we perform edge contractions if necessary so that there exist at most one vertex, say $v_+$ with $1={\rm indeg}(v_+)<{\rm outdeg}(v_+)$ and at most one vertex, say $v_-$ with ${\rm indeg}(v_-)>{\rm outdeg}(v_-)=1$. See for example Figure \ref{proof2} (e). Note that these replacements never decrease the number of connected components. Therefore it is sufficient to show that $S(D')$ is connected after these replacements.

%

%

%

Now suppose that there is just one sink of $G'$. Let $P$ be any point of $S(D')$. We start from $P$ along the flow of edge orientations of $S(D')$. If we come across the vertex $v_+$ then we choose the leftmost way. Namely we turn to the right at $v_+$. Then we see that as we turn around the origin of ${\mathbb R}^2$ we move to left and we finally reach to the sink. Thus we have that $S(D')$ is arcwise connected. Similarly if there are no sinks of $G'$ then starting from any point of $S(D')$ we finally reach to the outermost circle turning around the origin. Thus $S(D')$ is arcwise connected. Suppose that there are at most one source of $G'$. Then we see by going against the flow that $S(D')$ is arcwise connected.
Therefore it is sufficient to consider the case that there are at least two sinks and two sources of $G'$. Then by the definition of $G'$ we have that $G'$ has no vertices of degree one. 

\begin{figure}[htbp]
      \begin{center}
\scalebox{0.6}{\includegraphics*{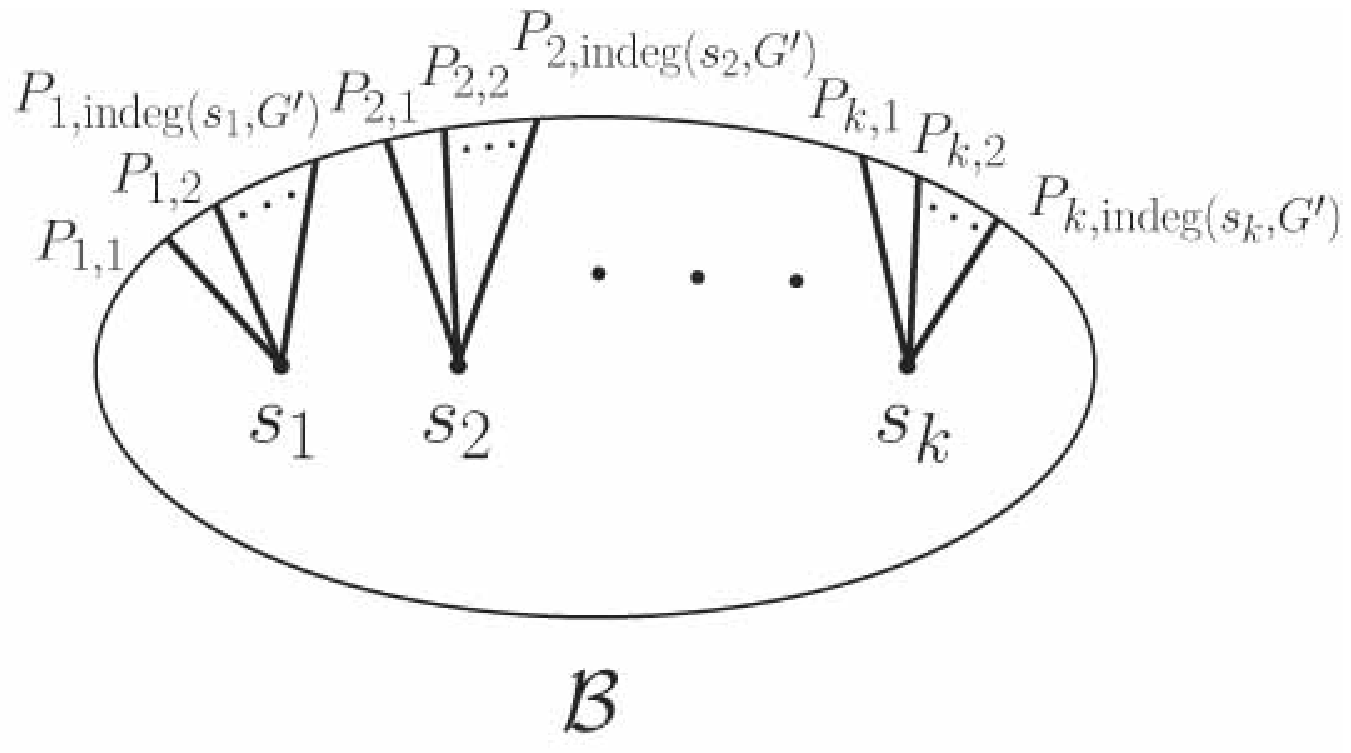}}
      \end{center}
   \caption{}
  \label{proof4}
\end{figure}

Let ${\mathcal B}$ be a disk in $B-B'$ containing all sinks in its interior. Let $s_1,\cdots,s_k$ be the sinks and $P_{1,1}$, $\cdots$, $P_{1,{\rm indeg}(s_1,G')}$, $\cdots$,$P_{k,1}$, $\cdots$, $P_{k,{\rm indeg}(s_k,G')}$ the points of intersection of $S(D')$ and $\partial{\mathcal B}$ such that they appear in this order on $\partial{\mathcal B}$, $P_{1,1}$ is adjacent to $v_-$ and $P_{i,j}$ is adjacent to $s_i$ for each $i$ and $j$. 
See for example Figure \ref{proof1} and Figure \ref{proof4}.

%

We will deform $D'$ only on ${\mathcal B}$. We divide the points $P_{1,1}$, $\cdots$, $P_{1,{\rm indeg}(s_1,G')}$, $\cdots$,$P_{k,1}$, $\cdots$, $P_{k,{\rm indeg}(s_k,G')}$ into some sets of points ${\mathcal S}_1,\cdots,{\mathcal S}_\alpha$ such that for each $i$ the points in ${\mathcal S}_i$ are consecutive on $\partial{\mathcal B}$ and any two points in ${\mathcal S}_i$ can be connected by an arc in $S(D')$ outside ${\rm int}{\mathcal B}$. We may suppose without loss of generality that $P_{1,1}\in{\mathcal S}_1$ and for each $i$ there is a pair of consecutive points on $\partial{\mathcal B}$ such that one is contained in ${\mathcal S}_i$ and the other is contained in ${\mathcal S}_{i+1}$ where we consider $\alpha+1=1$. 
We will show that ${\mathcal S}_i$ contains two or more points possibly except $i=1$. To see this we start from $P_{i,j}\neq P_{1,1}$ and trace $S(D')$ against the flow. Then we reach to $v_+$ or a source. Then we choose an edge that is next to the edge where we come and along the flow we trace $S(D')$. If we come across $v_+$ then we choose the leftmost way. 
Then we must reach to $P_{i,j+1}$ or $P_{i,j-1}$ where we consider $P_{i,0}=P_{i-1,{\rm indeg}(s_{i-1},G')}$, $P_{i,{\rm indeg}(s_i,G')+1}=P_{i+1,1}$ and $P_{k+1,1}=P_{1,1}$.

Let ${\mathcal B}={\mathcal B}_1\supset{\mathcal B}_2\supset\cdots\supset{\mathcal B}_k$ be a sequence of disks such that their boundaries $\partial{\mathcal B}_1,\partial{\mathcal B}_2,\cdots,\partial{\mathcal B}_k$ forms concentric circles in ${\mathcal B}$.
Let ${\mathcal A}_i={\mathcal B}_{i}-{\rm int}{\mathcal B}_{i+1}$ be the annulus.

\begin{figure}[htbp]
      \begin{center}
\scalebox{0.6}{\includegraphics*{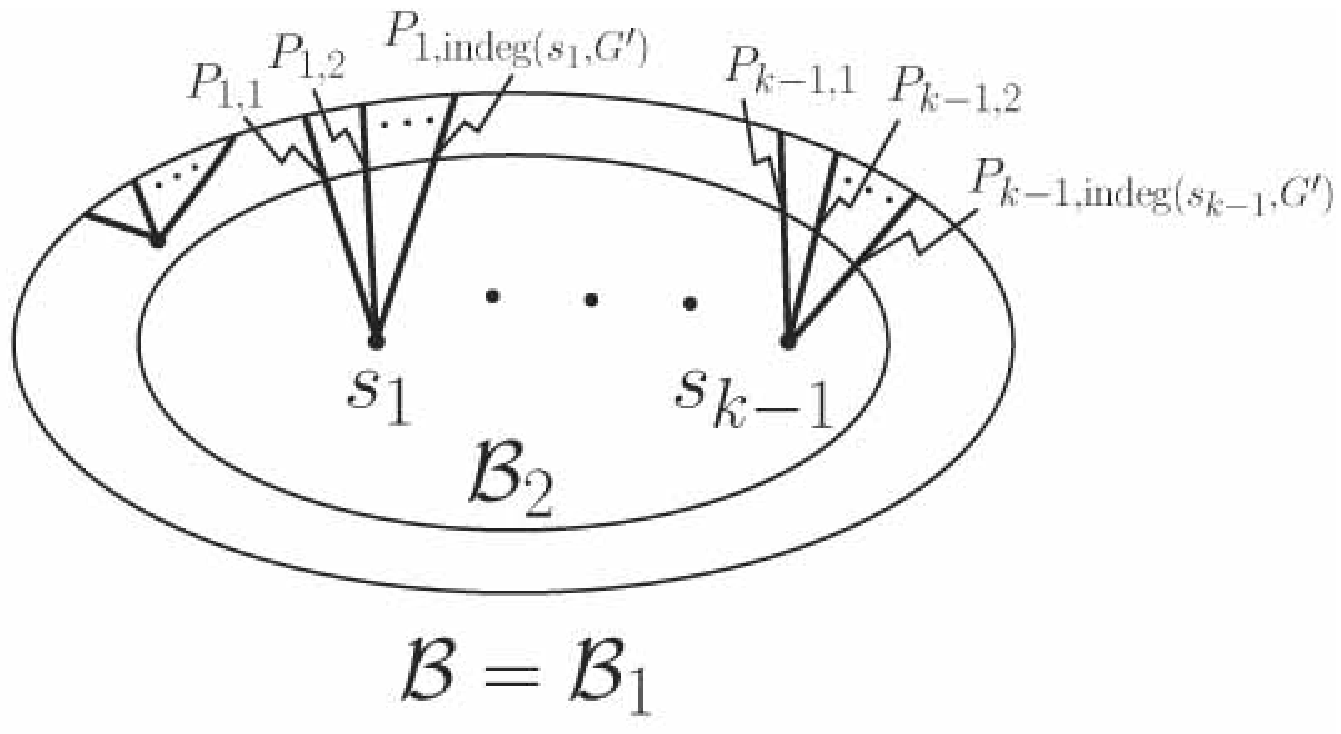}}
      \end{center}
   \caption{}
  \label{proof5}
\end{figure}

Suppose that the set $\{P_{1,1},\cdots,P_{1,{\rm indeg}(s_1,G')}\}$ is contained in ${\mathcal S}_1\cup\cdots\cup{\mathcal S}_i$ but not contained in ${\mathcal S}_1\cup\cdots\cup{\mathcal S}_{i-1}$.
First suppose that the set $\{P_{1,1},\cdots,P_{1,{\rm indeg}(s_1,G')}\}$ is a proper subset of ${\mathcal S}_1\cup\cdots\cup{\mathcal S}_i$.
Then we leave $s_1$ in ${\mathcal A}_1$ and rename the sinks $s_2,\cdots,s_k$ $s_1,\cdots,s_{k-1}$ and the points of intersection of $S(D')$ and $\partial{\mathcal B}_2$ as illustrated in Figure \ref{proof5}. Then we redivide the points $P_{1,1}$, $\cdots$, $P_{1,{\rm indeg}(s_1,G')}$, $\cdots$,$P_{k-1,1}$, $\cdots$, $P_{k-1,{\rm indeg}(s_{k-1},G')}$ into some sets of points, still denoted by ${\mathcal S}_1,\cdots,{\mathcal S}_\alpha$, such that for each $i$ the points in ${\mathcal S}_i$ are consecutive on $\partial{\mathcal B}_2$ and any two points in ${\mathcal S}_i$ can be connected by an arc in $S(D')$ outside ${\rm int}{\mathcal B}_2$. We may suppose without loss of generality that $P_{1,1}\in{\mathcal S}_1$ and for each $i$ there is a pair of consecutive points on $\partial{\mathcal B}_2$ such that one is contained in ${\mathcal S}_i$ and the other is contained in ${\mathcal S}_{i+1}$ where we consider $\alpha+1=1$.
Then by the construction we have that ${\mathcal S}_i$ contains two or more points possibly except $i=1$, or except $i=\alpha$.
If ${\mathcal A}_\alpha$ contains just one point then we reverse the cyclic order for the next step. Namely we rename again $s_1,s_2,\cdots,s_{k-1}$ $s_{k-1},s_{k-2},\cdots,s_{1}$ and rename ${\mathcal S}_1,{\mathcal S}_2,\cdots,{\mathcal S}_\alpha$ ${\mathcal S}_\alpha,{\mathcal S}_{\alpha-1},\cdots,{\mathcal S}_1$, and rename the points $P_{i,j}$ along the new cyclic order on $\partial{\mathcal B}_2$. 

\begin{figure}[htbp]
      \begin{center}
\scalebox{0.6}{\includegraphics*{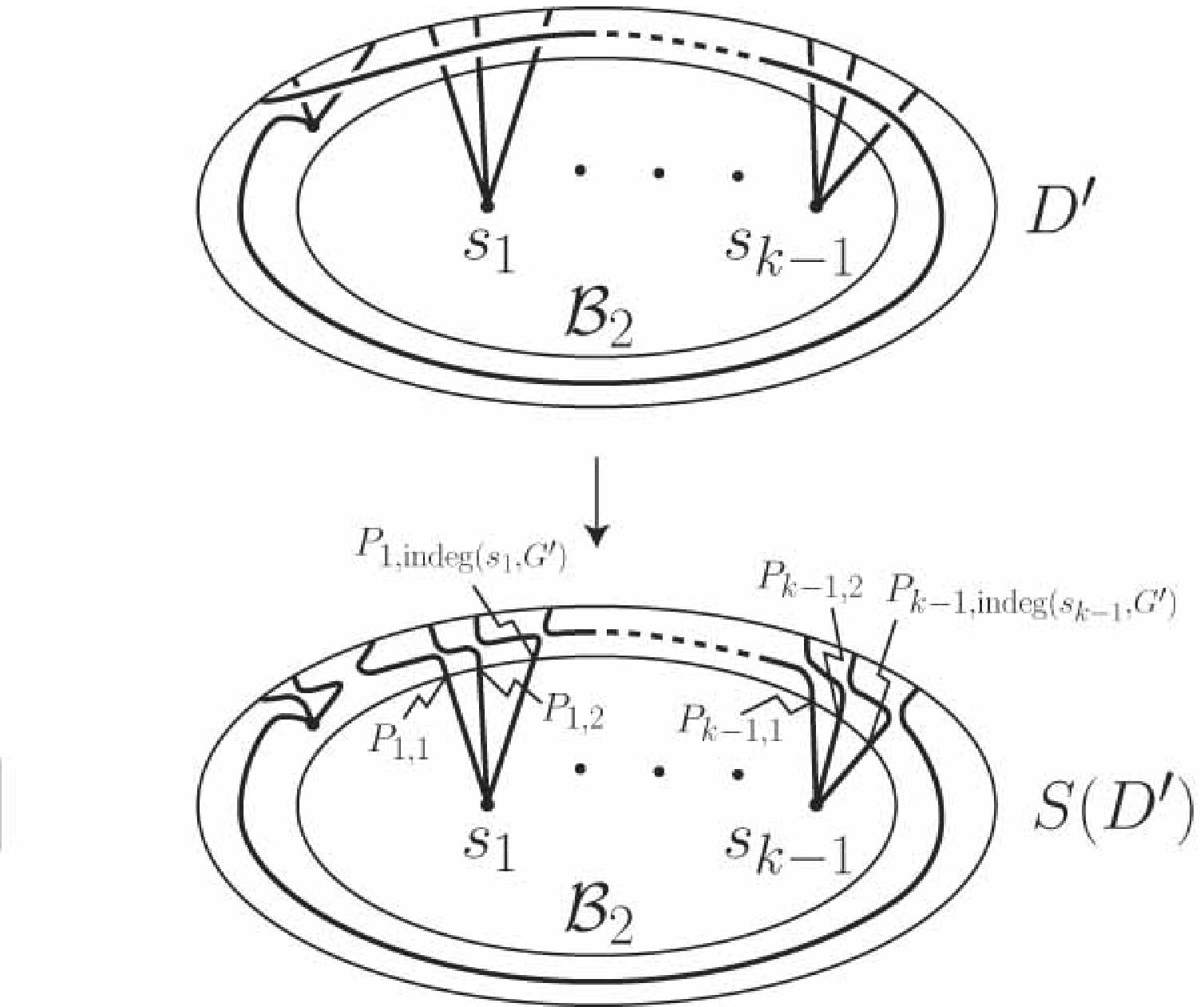}}
      \end{center}
   \caption{}
  \label{proof6}
\end{figure}

Next suppose that the set $\{P_{1,1},\cdots,P_{1,{\rm indeg}(s_1,G')}\}$ is equal to the set ${\mathcal S}_1\cup\cdots\cup{\mathcal S}_i$. 
Then we deform $D'$ as illustrated in Figure \ref{proof6} and consider $S(D')$. 
Note that new $P_{1,1}$ and new $P_{k-1,{\rm indeg}(s_{k-1},G')}$ can be connected by an arc in $S(D')$ outside ${\rm int}{\mathcal B}_2$.
Therefore we have that each new ${\mathcal S}_i$ contains at least two points. 
Next we deform $D'$ inside ${\mathcal B}_2$ and leave new $s_1$ in ${\mathcal A}_2$ in a similar way.
We continue this deformation and finally have the desired $S(D')$. 

\begin{figure}[htbp]
      \begin{center}
\scalebox{0.6}{\includegraphics*{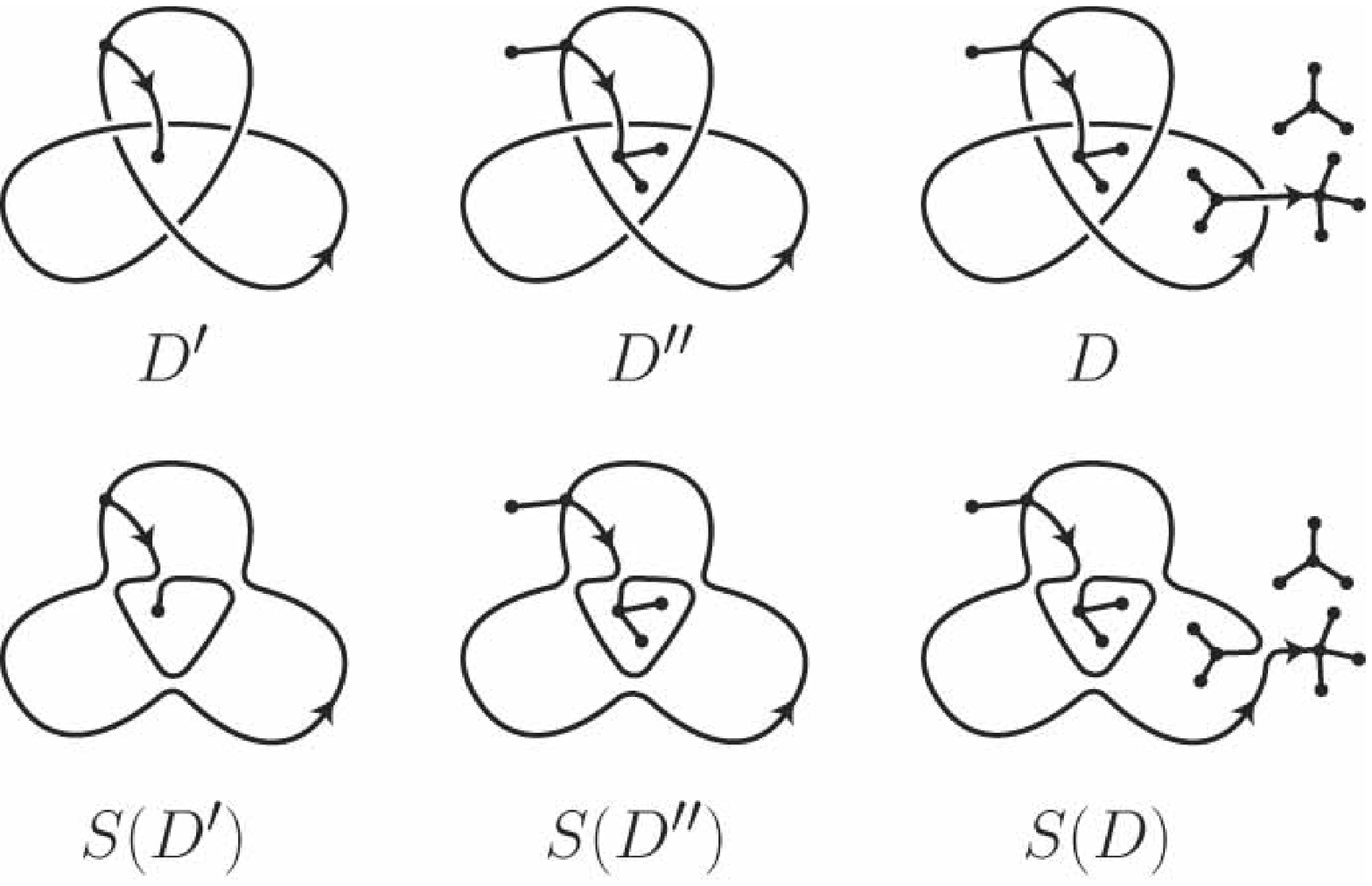}}
      \end{center}
   \caption{}
  \label{proof7}
\end{figure}

Now we return to the whole graph $G$. Let $G''$ be the maximal subgraph of $G$ that contains $G'$ and $\mu(G'')=\mu(G')$. Let $T_1,\cdots,T_n$ be the tree components of $G$ that are disjoint from $G'$. Then $G=G''\cup T_1\cup\cdots\cup T_n$.
Let $f''$ be the restriction of the spatial embedding $f:G\to{\mathbb S}^3$ to $G''$. Let $D''$ be a diagram of $f''$ whose subdiagram for $f'$ is $D'$ and has no more crossings than $D'$. Then we have that $S(D'')$ and $S(D')$ have the same homotopy type. In particular $S(D'')$ is connected. Let $m={\rm min}(\beta_1(S(D'')),n)$. Let $Q_1,\cdots,Q_m$ be points on $S(D'')$ other than the vertices such that $S(D'')-\{Q_1,\cdots,Q_m\}$ is still connected. We may suppose that these points are away from the neighbourhoods of the crossings of $D''$ where the smoothings are performed. Let $D$ be a diagram of $f$ whose subdiagram for $f''$ is $D''$ such that the crossings of $D$ other than that of $D''$ are exactly the points $Q_1,\cdots,Q_m$ where the crossing $Q_i$ is between an edge of $G'$ and an edge of $T_i$. 
Then we see that $\mu(S(D))=1+n-m$. See for example Figure \ref{proof7}. 
Note that we have the following equality.
\[
\chi(G)=\chi(G'')+n=\chi(S(D''))+n=\mu(S(D''))-\beta_1(S(D''))+n=1-\beta_1(S(D''))+n.
\]
Therefore if $m={\rm min}(\beta_1(S(D'')),n)=n$ then we have $\chi(G)\leq1$ and $\mu(S(D))=1+n-m=1$ as desired.
If $m={\rm min}(\beta_1(S(D'')),n)=\beta_1(S(D''))$ then we have $\chi(G)\geq1$ and $\mu(S(D))=1+n-\beta_1(S(D''))=\mu(S(D''))-\beta_1(S(D''))+n=\chi(S(D''))+n=\chi(G'')+n=\chi(G)$ as desired. This completes the proof. 
$\Box$

%

%

%

%

%

\vskip 3mm

\begin{figure}[htbp]
      \begin{center}
\scalebox{0.6}{\includegraphics*{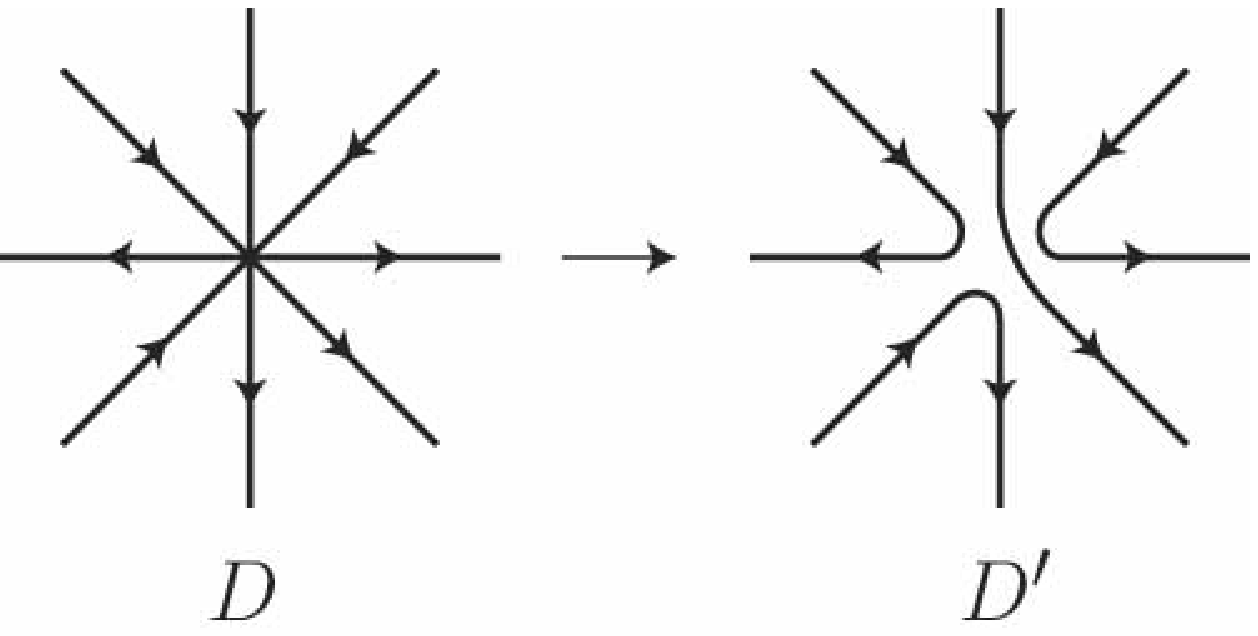}}
      \end{center}
   \caption{}
  \label{proof8}
\end{figure}

\begin{figure}[htbp]
      \begin{center}
\scalebox{0.6}{\includegraphics*{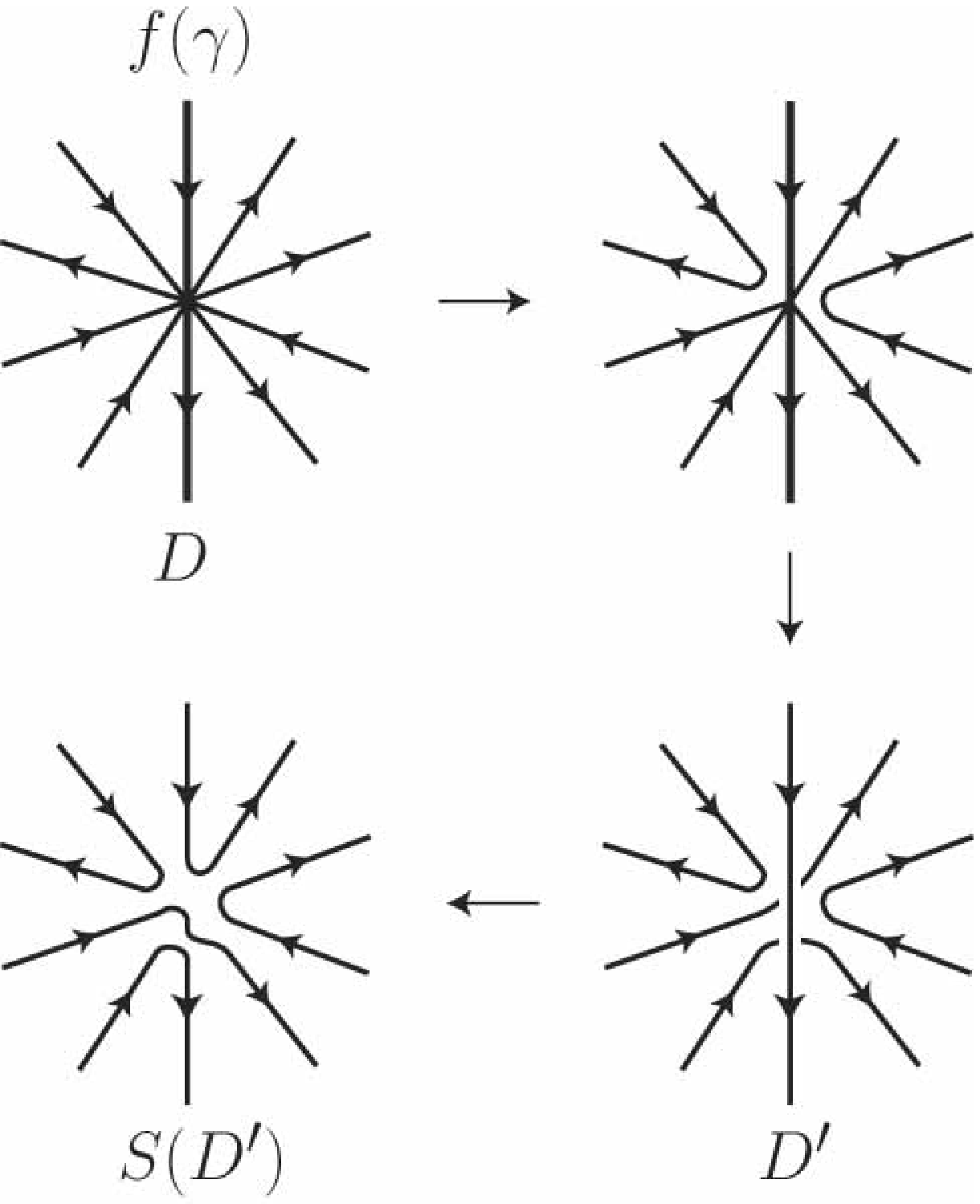}}
      \end{center}
   \caption{}
  \label{proof9}
\end{figure}

\noindent{\bf Proof of Theorem \ref{main-theorem2} (2).} Let $\gamma$ be an oriented cycle of $G$. Let $f:G\to {\mathbb S}^3$ be a spatial embedding of $G$ such that the braid index $b(f(\gamma))$ of the knot $f(\gamma)$ is greater than or equal to $n-\chi(G)$.
Let $D$ be any diagram of $f$. It is sufficient to show that $\mu(S(D))\geq n$. We replace each neighbourhood of a vertex $v$ of $D$ to ${\rm indeg}(v,G)$  oriented arcs as follows. For a vertex $v$ that is not on $\gamma$ we replace it to mutually disjoint oriented arcs. See for example Figure \ref{proof8}.
Let $v$ be a vertex of $G$ that is on $\gamma$. Let $N$ be a small neighbourhood of $v$ on $D$. Suppose that there is a pair of edges not contained in $\gamma$, say $e_i$ and $e_o$, such that the head of $e_i$ is $v$, the tail of $e_o$ is $v$ and they are next to each other in $N$. Then we take them away from $v$ and connect them. We do this for all such pairs. Then we have the situation that all edges in $N$ not on $f(\gamma)$ go from the right of $f(\gamma)$ to the left of $f(\gamma)$ or from the left of $f(\gamma)$ to the right of $f(\gamma)$. Then we split off them and let $f(\gamma)$ goes over them. Let $D'$ be the result of these replacements. 
See for example Figure \ref{proof9}. Then we have that $D'$ is a diagram of some oriented link, say $L$. Since $L$ contains a knot $f(\gamma)$ we have that the braid index $b(L)$ of $L$ is greater than or equal to $n-\chi(G)$. By the result in \cite{Yamada} we have that $\mu(S(D'))\geq b(L)$. Therefore we have $\mu(S(D'))\geq n-\chi(G)$. Note that the homotopy type of $S(D)$ is obtained from $S(D')$ by adding $\sum_{v}({\rm indeg}(v,G)-1)$ edges where the summation is taken over all vertices $v$ of $G$. Therefore we have that 
$\mu(S(D))\geq\mu(S(D'))-\sum_{v}({\rm indeg}(v,G)-1)$. By the handshaking lemma and by the assumption that $G$ is circulating we have that ${\displaystyle \sum_{v}({\rm indeg}(v,G)-1)=-\chi(G)}$. 
Thus we have $\mu(S(D))\geq n$ as desired. $\Box$

%

%

%

The following example shows that even for the circulating graphs the difference $b(f)-s(f)$ depends on the spatial embedding $f$.

\vskip 3mm

\begin{figure}[htbp]
      \begin{center}
\scalebox{0.4}{\includegraphics*{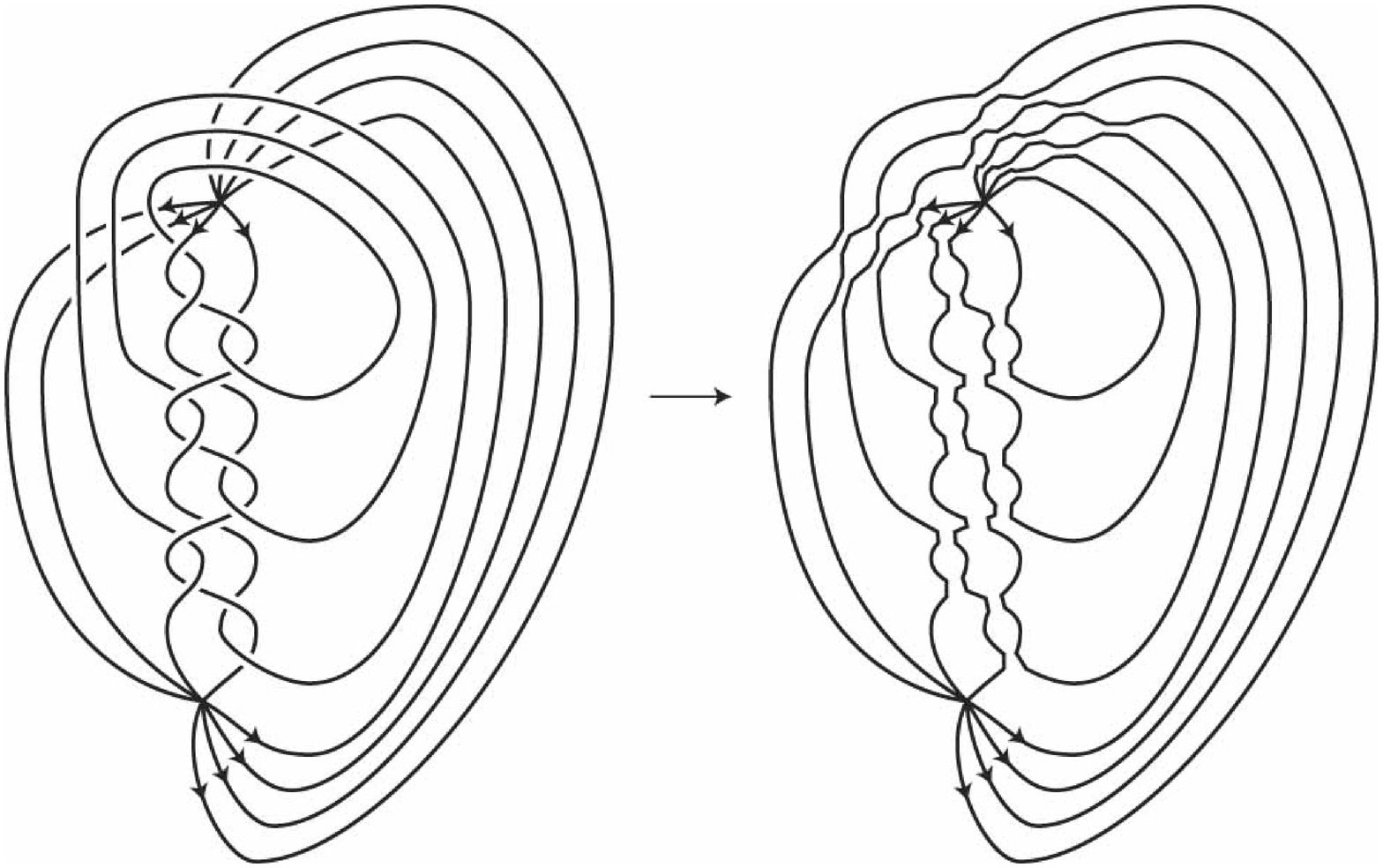}}
      \end{center}
   \caption{}
  \label{example1}
\end{figure}

\begin{figure}[htbp]
      \begin{center}
\scalebox{0.6}{\includegraphics*{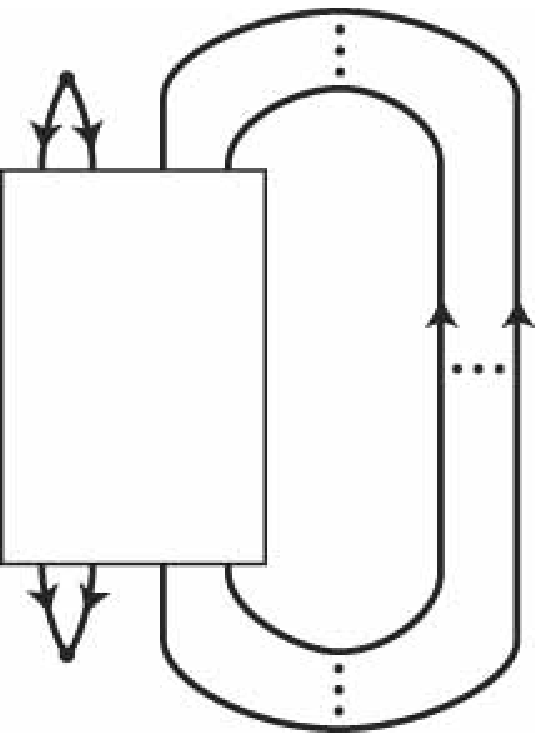}}
      \end{center}
   \caption{}
  \label{example2}
\end{figure}

\begin{Example}\label{example}
{\rm
Let $G$ be a circulating graph on two vertices and eight edges joining them. 
Let $f:G\to {\mathbb S}^3$ be a trivial embedding of $G$. Then we have that $b(f)=4$ and $s(f)=1$.
Let $g:G\to {\mathbb S}^3$ be a spatial embedding of $G$ illustrated in Figure \ref{example1}. Note that $g(G)$ contains a knot $K$ that is a connected sum of three figure eight knots. Then we have ${\rm bridge}(K)=4$. Suppose that $K$ is a braid presentation as its edge orientations. Then we may suppose that $K$ is as illustrated in Figure \ref{example2} where the box represents some $n$ braid. Then we have that ${\rm bridge}(K)\leq n-1$. Therefore we have that $n\geq5$. Since $G$ has two more oriented cycles other than $K$ we have that $b(g)\geq7$. Since $g$ is a braid presentation with $\tilde{b}(g)=7$ we have $b(g)=7$. However we have that $s(g)=1$ as illustrated in Figure \ref{example1}.
}
\end{Example}

%

%

%

\section*{Acknowledgments}
The authors are grateful to Professor Shin'ichi Suzuki for his constant guidance and encouragement. The authors are also grateful to Dr. Ryuzo Torii for his helpful comments.

{\normalsize
}

\end{document}